\documentclass[final,5p,times]{elsarticle}

\usepackage{amssymb}
\usepackage{amsthm}
\usepackage{amsmath}
\usepackage{lipsum} 
\usepackage[ruled,vlined]{algorithm2e}
\usepackage{multirow,tabularx}
\usepackage{caption}
\usepackage{subcaption}
\usepackage{xcolor}
\usepackage{natbib}
\usepackage{hyperref}

\newdefinition{defn}{Definition}
\newdefinition{assn}{Assumption}
\newtheorem{thm}{Theorem}
\newtheorem{prop}{Proposition}

\newdefinition{rmk}{Remark}
\newproof{prf}{Proof}

\usepackage{mathtools}

\newcommand{\tb}[1]{{\color{black}#1}}

\begin{document}

\begin{frontmatter}



\title{Robust Confidence Bands for Stochastic Processes Using Simulation}


\author{Timothy C. Y. Chan}

\author{Jangwon Park\corref{cor1}}
\cortext[cor1]{Corresponding author at: Department of Mechanical and Industrial Engineering, University of Toronto, Toronto, Ontario, Canada M5S 3G8.}
\ead{jangwon.park@mail.utoronto.ca}

\author{Vahid Sarhangian}

\address{Department of Mechanical and Industrial Engineering, University of Toronto, Canada}

\begin{abstract}
We propose a robust optimization approach for constructing confidence bands for stochastic processes using a finite number of simulated sample paths. Our methodology addresses optimization bias within the constraints, avoiding overly narrow confidence bands of existing methods. In our first case study, we show that our approach achieves the desired coverage rates with an order-of-magnitude fewer sample paths than the state-of-the-art baseline approach. In our second case study, we illustrate how our approach can validate stochastic simulation models.
\end{abstract}




\end{frontmatter}




\section{Introduction} \label{sec:intro}
Stochastic simulation \cite{asmussen2007stochastic,nelson2021foundations} is a primary tool for performance evaluation of stochastic dynamical systems, in particular under counterfactual or ``what-if" scenarios.  In many applications, simulation outputs are sample paths of stochastic processes realized over a finite horizon. For instance, these sample paths may correspond to hospital occupancy levels \cite{chan2021optimizing, helmdelta}, the price of financial products \cite{lamberton2011introduction}, or the number of infected patients in an infectious disease model \cite{cramer2022evaluation}. A natural way to validate a simulation model in this setting is to construct a \textit{confidence band} over the sample paths at a specified level of coverage rate, and check whether the historical paths from the actual system are ``covered'' by the confidence band (see Section \ref{sec:prelim} for the formal definition). 

In addition to validation, confidence bands can also be used to quantify the uncertainty in realizations of a stochastic process over a finite horizon. This has received much attention in the context of estimating impulse response functions in vector autoregressive (VAR) models in the Economics literature \citep[e.g.,][]{simsErrorBandsImpulse1999, staszewskaRepresentingUncertaintyResponse2007, jordaSimultaneousConfidenceRegions2009, staszewska-bystrovaConstructingNarrowestPathwise2013, lutkepohlComparisonMethodsConstructing2015, schusslerConstructingMinimumwidthConfidence2016}. Confidence bands on impulse response functions are commonly estimated to examine the effects of shocks to a system over time. However, the majority of approaches are either heuristics or asymptotic methods, which may provide an overly wide confidence band when the number of observations (sample paths) is small. Moreover, these approaches are not universally applicable given their specific application to VAR. The best known approach that is most relevant to our work is \cite{schusslerConstructingMinimumwidthConfidence2016}, who propose a mixed-integer program (MIP) to construct minimum-width confidence bands. As we show in this work, however, this approach may produce a biased confidence band, i.e., one with a smaller coverage rate than desired. Although generating more sample paths reduces this bias, the number of sample paths required to achieve a sufficiently small bias may result in a large MIP that is not solvable to optimality in practical time, as noted in \cite{schusslerConstructingMinimumwidthConfidence2016}.

Related methods for specific stochastic processes include \cite{kendallConfidenceBandsBrownian2007}, who formulate an optimization problem based on local time arguments for constructing confidence bands on Brownian motion and perturbed Brownian paths. While the approach is relevant when considering Brownian approximations, it is not applicable to general simulation output that may not be well-approximated by Brownian motion. On the other hand, \cite{korpelaMultivariateConfidenceIntervals2017} and \cite{bergMinimumWidthConfidenceBands2017} propose various heuristics and an MIP approach, respectively, which are applicable for stochastic processes such as stock market data, temperature data, and medical data (e.g., heartbeat). However, these works pursue an alternative definition of coverage, which is arguably less practical and less amenable to theoretical analysis; see also the discussion in Section \ref{sec:prelim}. Therefore, we do not consider it here. Furthermore, they do not directly address the potential optimization bias that may lead to an overly narrow confidence band.

Also related is the literature on conformal prediction \citep[e.g.,][]{leiConformalPredictionApproach2015, romanoConformalizedQuantileRegression2019, stankeviciuteConformalTimeseriesForecasting2021, zaffranAdaptiveConformalPredictions2022, xuConformalPredictionTime2023, auerConformalPredictionTime2023} and prediction intervals for metamodels \citep[e.g.,][]{lamPredictionIntervalsSimulation2022}. In this context, the aim is to construct prediction intervals that capture both the intrinsic (aleatoric) and extrinsic (epistemic) uncertainty, arising from model mis-specification, limited data, and the inherent stochasticity of the model. As such, these methods typically involve training and validation phases to compute the residual errors of predictors to calibrate the width of the intervals. Our work is different because we focus only on the intrinsic uncertainty and more data can be generated at will (i.e., by simulation). Thus, we cast the task of constructing intervals as an empirical constrained optimization problem where we minimize the width of the intervals directly, subject to a desired coverage rate over a set of generated sample paths. In \cite{lamPredictionIntervalsSimulation2022}, the authors  consider a similar problem for simulation metamodeling, whose solution is the prediction interval for the response surface $\mathbb{E}[Y(x)]$, where $x$ is an input parameter (e.g., time) and $Y(x)$ is a random output from simulation (e.g., number of customers). The main difference in our work is the ability to solve the empirical constrained optimization problem directly, rather than an approximation of it, due to our specific focus on (discrete-time) stochastic processes and the intrinsic uncertainty. Additionally, while \cite{lamPredictionIntervalsSimulation2022} adds a positive constant to the coverage rate to avoid overly narrow intervals, we add a positive constant to the total width of the intervals instead, which allows more flexibility since one can increase the width without affecting the coverage rate, but not vice versa. 


In this paper, we propose a new approach for constructing confidence bands for general discrete-time stochastic processes using a finite number of simulated sample paths. Unlike existing approaches in the literature, our methodology is widely applicable and directly addresses optimization bias through a robust optimization approach. It is tractable, being only slightly more complex than the state-of-the-art baseline approach, and easy to use, as it employs standard techniques. Our approach is also applicable to continuous-time processes after appropriately discretizing time. \tb{While our main theoretical results focus on continuous-valued processes, the method remains applicable in practice to discrete-valued processes with large support.} In our first case study, we demonstrate that our confidence bands achieve the desired coverage rates with an order-of-magnitude fewer sample paths than the best known approach for the same problem. In the second case study, we illustrate how our method can help validate stochastic simulation models. The data and code for our numerical experiments are available on \href{https://github.com/parkjan4/RobustConfidenceBands}{GitHub}\footnote{https://github.com/parkjan4/RobustConfidenceBands}.


\section{Preliminaries} \label{sec:prelim}
Let $\{X_t; t\geq 0\}$ denote a (discrete-time) stochastic process and let $X \equiv (X_1,\ldots,X_H)$ denote the finite-dimensional random vector corresponding to the value of the process at times $\{1,\ldots,H\}$, taking values in $\mathbb{R}^H$. The path $X$ could, for example, represent the output of a simulation model that we are interested in validating, or an impulse response function. \tb{Throughout this paper, we make the following assumptions on $X_t$.}
{\color{black}
\begin{assn} \label{assn:bounds}
    There exist finite $L_t$ and $U_t$ such that $P(L_t \leq X_t \leq U_t)=1$ for all $t\in \{1,\ldots,H\}$.
\end{assn}
\begin{assn} \label{assn:continuous_cdf}
    $X_t$ has a continuous cumulative distribution function (CDF) for all $t\in \{1,\ldots,H\}$.
\end{assn}
}
\tb{Assumption~\ref{assn:bounds} is not practically restrictive, as in most applications we can find natural upper and lower bounds. For instance, for hospital occupancy, $ U_t $ can be set to the hospital's capacity and $ L_t = 0 $ for all $ t $. In other contexts where sample paths are technically unbounded, e.g.,  financial prices, one can set sufficiently large but practical upper (lower) bounds to ensure Assumption~\ref{assn:bounds} holds. Assumption~\ref{assn:continuous_cdf} precludes discrete random variables. While this assumption is required for our main results, we show in Section~\ref{ssec:discrete_case} that the proposed methods are also applicable to practical settings with discrete random variables that take values in a sufficiently large set.}

Denote a confidence band by $(l,u)$ where $l,u\in \mathbb{R}^H$. We formally define ``coverage'' as follows.
\begin{defn}[Coverage] \label{defn:coverage}
    We say $x \equiv (x_1,\ldots,x_H) \in \mathbb{R}^H$ is \textit{covered} by $(l,u)$ if $l_t \leq x_t \leq u_t$ for all $t=1,\ldots,H$.
\end{defn}
Definition \ref{defn:coverage} is a natural definition for coverage. We do not pursue the alternative definition in \cite{korpelaMultivariateConfidenceIntervals2017} and \cite{bergMinimumWidthConfidenceBands2017}, which states that a sample path is covered if it does not lie outside the confidence band by more than $l$ times (where $l$ is user-specified). Although more general, it is less practical since it is not possible to know at \textit{which} time steps a given sample path will lie inside the confidence band. Moreover, computing its coverage rate requires considering all possible combinations of $s\in \{0,\ldots,l\}$ time steps that lie outside the band, which is less conducive to theoretical analysis, especially when coverage at different time steps is correlated.

Let $C_t = \{l_t \leq X_t \leq u_t\}$ be the event that $X$ is covered by the given confidence band $(l,u)$ at time $t$. We make the following assumption on $C_t$. 
\begin{assn} \label{assn:coverage_uncertainty}
    For any $t'$ and $t$ such that $t'>t$, $P(C_{t'}|C_t) < 1$.
\end{assn}
Assumption \ref{assn:coverage_uncertainty} is not restrictive. It states that there is uncertainty in coverage over time so that if $X$ is covered at time $t$, this does not guarantee coverage at a future time. In the rest of the paper, for ease of exposition, we use the shorthand notation $P(l,u) = P(C_1, \ldots, C_H)$ to denote the coverage rate of $(l,u)$. We omit its dependence on $H$, as it should be clear from context.

When $H=1$, the confidence band is simply the interval enclosed between the $(1-\alpha/2)$- and $(\alpha/2)$-quantiles, where $1-\alpha \in [0,1]$ is the specified coverage rate. When $H>1$, one may consider the ``naive confidence band'' by extending the $H=1$ case as follows: estimate the $(1-\alpha/2)$- and $(\alpha/2)$-quantiles at each time step, and connect the upper and lower points of the adjacent intervals. Let $(\Bar{l}, \Bar{u})$ denote the naive confidence band using exact quantiles at all time points. Although intuitive and often used in applied work, this solution provides a smaller coverage rate than desired, as formalized below.

\begin{prop} \label{prop:naiveCB}
    Let $1-\alpha$ be the desired coverage rate. Under Assumptions \ref{assn:continuous_cdf} and \ref{assn:coverage_uncertainty}, the naive confidence band under-covers, i.e., $P(\Bar{l}, \Bar{u}) < 1-\alpha$.
\end{prop}
\begin{proof}
    The proof is by induction on $H$, the length of the horizon. Let $C_t$ be the event that $X$ is covered by $(\Bar{l}, \Bar{u})$ at time $t$. Suppose $H=2$. We observe that
    \begin{align*}
        P(\Bar{l}, \Bar{u}) = P(C_1,C_2) = P(C_2|C_1)P(C_1) = P(C_2|C_1)(1-\alpha) < 1-\alpha,
    \end{align*}
    where \tb{$P(C_1)=1-\alpha$ is implied by Assumption \ref{assn:continuous_cdf}}, and  the last inequality holds by Assumption \ref{assn:coverage_uncertainty}. Now, suppose the same holds for $H=n$. Then for $H=n+1$, we observe that
    \begin{align*}
        P(\Bar{l}, \Bar{u}) = P(C_1,\ldots,C_{n+1}) &= P(C_{n+1}|C_1,\ldots,C_n)P(C_1,\ldots,C_n) \\
        &< P(C_{n+1}|C_1,\ldots,C_n)(1-\alpha) \\
        &< 1-\alpha,
    \end{align*}
    where the first inequality is by the induction hypothesis and the second inequality is again by Assumption \ref{assn:coverage_uncertainty}.
\end{proof}

\section{Baseline: nominal MIP} \label{sec:baseline}
We first present the baseline model in \cite{schusslerConstructingMinimumwidthConfidence2016}, which we refer to as the \textit{nominal problem}. Suppose we have generated $n$ iid sample paths. Let $x^i \in \mathbb{R}^H$ be the $i$th sample path and denote by $1-\alpha$ the desired coverage rate. Let $\delta(l,u,x^i)$ denote an indicator function that evaluates to 1 if and only if $x^i$ is covered by the confidence band $(l,u)$, where $l = (l_1, \ldots, l_H)$ and $u=(u_1, \ldots, u_H)$. Denote by $q^u_t$ and $q^l_t$ the $(1-\alpha/2)$- and $(\alpha/2)$-quantile estimates at time $t$ based on the $n$ sample paths. The nominal problem finds a minimum-width confidence band such that it covers at least $(1-\alpha)\times 100\%$ of the sample paths:

\begin{equation}\label{nominal_mip} \tag{NP}
    \tb{\begin{split}
        w^\star_n = \min_{(l,u) \in \mathcal{X}} \quad & \sum_{t=1}^H(u_t - l_t) \\
        \mbox{s.t.} \quad & \frac{1}{n}\sum_{i=1}^{n}\delta(l,u,x^i) \geq 1-\alpha,
    \end{split}}
\end{equation}
where \tb{$\mathcal{X} = \{(l,u): u_t \geq q^u_t, l_t \leq q^l_t, \forall t\}$}. This problem can be reformulated as a MIP using $n$ binary variables and an appropriate large constant as shown in \cite{schusslerConstructingMinimumwidthConfidence2016}. Hereafter, we refer to an optimal solution of \eqref{nominal_mip} as the \textit{nominal confidence band}. We focus on the case where $n(1-\alpha)$ is an integer. Otherwise, the constraint $(1/n)\sum_{i=1}^{n}\delta_i \geq 1-\alpha$ holds strictly, and one may consider the smallest $\beta \in (\alpha, 1]$ such that $n(1-\beta)$ is an integer.

\tb{In the rest of this section, we analyze the properties of the nominal confidence band. We first observe that \eqref{nominal_mip} can be viewed as a sample average approximation to the following optimization problem with the ``true'' probabilistic constraint:}
\tb{\begin{equation}\label{true_problem} \tag{$\star$}
    \begin{split}
        w^\star = \min_{(l,u) \in \mathcal{X}} \quad & \sum_{t=1}^H(u_t - l_t) \\
        \mbox{s.t.} \quad & P(l, u) \geq 1-\alpha,
    \end{split}
\end{equation}
where $P(l,u)$ is a shorthand notation for $P(l_t \leq X_t \leq u_t, \forall t)$. Let $(l^\star, u^\star)$ be an optimal solution to \eqref{true_problem}, and let $w^\star$ be its total width. We note that an optimal solution to \eqref{true_problem} exists, since under Assumptions \ref{assn:bounds} and \ref{assn:continuous_cdf}, $P(l,u)$ is continuous in $l$ and $u$, and the feasible set of \eqref{true_problem} is non-empty and compact. We now state the following result about the nominal confidence band.}
\tb{\begin{thm} \label{thm:undercoverage}
      $P(w^\star_n \leq w^\star) > \frac{1}{2}$ for all $n\in \mathbb{Z}_+$. Moreover, for any $(l,u) \in \mathcal{X}$ with $\sum_{t=1}^H (u_t-l_t) \leq w^\star$, $P(l, u) \leq 1-\alpha$.
\end{thm}}
\tb{\begin{proof}
    To prove the first part, we use Lemma 1 from \cite{luedtkeSampleApproximationApproach2008}, who study sample average approximations of optimization problems with probabilistic constraints. For $\epsilon,\alpha \in [0,1]$, a random vector $\xi$, a given function $G$, and a deterministic set $\mathcal{D}$, let
    \begin{align*}
        \mathcal{D}_\epsilon &= \left\{x\in \mathcal{D}: P(G(x,\xi)\leq 0) \geq 1-\epsilon \right\}, \\ 
        \mathcal{D}^n_\alpha &= \left\{x\in \mathcal{D}: \frac{1}{n}\sum_{i=1}^n 1\left\{G(x,\xi^i) \leq 0 \right\} \geq 1-\alpha \right\},
    \end{align*}
    where $1\{\cdot\}$ is the indicator function. For a given objective function $f(x)$, let
    \begin{align*}
        z^\star_\epsilon &= \min\{f(x): x\in \mathcal{D}_\epsilon\}, \\
        \hat{z}^n_\alpha &= \min\{f(x): x\in \mathcal{D}^n_\alpha \}.
    \end{align*} 
    Lemma 1 of \cite{luedtkeSampleApproximationApproach2008} establishes that
    \begin{align*}
        \tb{P(\hat{z}^n_\alpha \leq z^\star_\epsilon) \geq \sum_{i=0}^{n\alpha} \binom{n}{i}\epsilon^i (1-\epsilon)^{n-i}.}
    \end{align*} 
    To use this result, we must express $P(l,u) = P(l_t \leq X_t \leq u_t, \forall t)$ as $P(G(l,u,X) \leq 0)$ for some function $G$. To this end, let $f_t(l,X) \equiv l_t - X_t$ and $g_t(u,X) \equiv X_t - u_t$ for $t\in \{1,\ldots, H\}$. Now define $G(l,u,X) \equiv \max_{1\leq t\leq H}\{f_t(l,X), g_t(u,X)\}$. Then it is straightforward to see that $G(l,u,X) \leq 0$ if and only if $f_t(l,X) \leq 0$ and $g_t(u,X) \leq 0$ for all $t$, or equivalently $l_t \leq X_t \leq u_t$ for all $t$.
    
    Specializing Lemma 1 of \cite{luedtkeSampleApproximationApproach2008} for the case of $\epsilon=\alpha$, we have
    \begin{align*}
        P(w^\star_n \leq w^\star) \geq F(n\alpha),
    \end{align*}
    where $F(n\alpha) \equiv \sum_{k=1}^{n\alpha} \binom{n}{k} \alpha^k (1-\alpha)^{n-k}$ is the CDF of the binomial distribution with $n$ trials and success probability $\alpha$. If $n\alpha$ is an integer, as we assumed earlier in Section \ref{sec:baseline}, it is well-known that $n\alpha$ is the (unique) strong median, i.e., $P(Y\leq n\alpha) > 1/2$ and $P(Y \geq n\alpha) > 1/2$ for $Y \sim \mbox{Binomial}(n,\alpha)$ \citep[Theorem~1]{kaasMeanMedianMode1980}. Thus, $P(w^\star_n \leq w^\star) \geq F(n\alpha) > 1/2$.
    
    The second part of the statement follows from the fact that any $(l,u) \in \mathcal{X}$ with $\sum_{t=1}^H (u_t-l_t) < w^\star$ must satisfy $P(l,u) < 1-\alpha$, since otherwise it would violate the optimality of $(l^\star, u^\star)$.
\end{proof}}
\tb{Theorem~\ref{thm:undercoverage} formalizes the under-coverage bias of the nominal confidence band: it is often overly narrow, leading to a coverage rate below $1 - \alpha$.} Intuitively, the nominal problem ``over-optimizes" with respect to the given samples and produces an optimistic (narrow) confidence band. This phenomenon is observed in other stochastic optimization contexts and referred to as \emph{optimization bias} \cite{homem2014monte} or the \emph{optimizer's curse} \cite{smith2006optimizer}. Although this bias diminishes as the number of sample paths increases, the resulting MIP to achieve a sufficiently small bias may be too large to solve to optimality in practical time. This motivates the development of a new approach to generate confidence bands that can achieve the desired coverage rates even with relatively small numbers of samples.

\tb{\subsection{The discrete case} \label{ssec:discrete_case}} 
\tb{In closing this section, we discuss the case when $X_t$ takes values only on a discrete set at each $t$. In such settings, the nominal confidence band's coverage rate may exceed $1-\alpha$ due to duplicate sample paths, since small decreases in the width of the band can significantly reduce coverage. Nonetheless, we show that under-coverage bias persists when the support of $X_t$ is reasonably large, as is common in many practical applications.}

\tb{Let $S_t$ denote the support of $X_t$, and let $\mathrm{supp}(X) = S_1 \times \cdots \times S_H$ denote the support of the random vector $X$, corresponding to all combinations of the discrete values at each $t$. For $s \in \mathrm{supp}(X)$, define \textit{multiplicity} $c(s) \equiv \#\{i: x^i_t=s_t, \forall t\}$, which counts the number of times $s$ was sampled, and let $c_{\mathrm{max}} \equiv \max_s c(s)$ be the maximum count. We next establish an upper bound on the nominal confidence band's coverage rate.
\begin{prop} \label{lem:bound}
    Let $(\hat{l}, \hat{u})$ denote an optimal solution to \eqref{nominal_mip}. Then $\frac{1}{n} \sum_{i=1}^n \delta(\hat{l}, \hat{u}, x^i) \leq 1 - \alpha + (c_{\mathrm{max}}/n)$.
\end{prop}
\begin{proof}
    For ease of notation, let $K = \frac{1}{n} \sum_{i=1}^n \delta(\hat{l}, \hat{u}, x^i)$. While $K > 1-\alpha$ may hold if $X_t$ takes on discrete values at each $t$, here we show that it cannot exceed $1-\alpha + (c_{\mathrm{max}}/n)$.

    Suppose $K > 1-\alpha + (c_{\mathrm{max}}/n)$. Pick an index $t$ such that for some $i$ satisfying $\hat{l}_t \leq x^i_t \leq \hat{u}_t$ for all $t$, either $\hat{l}_t=x^i_t$ or $\hat{u}_t=x^i_t$. Such $t$ exists by construction, since it contradicts minimal width of $(\hat{l}, \hat{u})$ otherwise. Consider shrinking $(\hat{l}, \hat{u})$ infinitesimally at $t$. Then its coverage rate is reduced by $c(x^i) / n$ but we still have 
    \begin{align*}
        K - \frac{c(x^i)}{n} > 1 - \alpha + \frac{c_{\mathrm{max}} - c(x^i)}{n} \geq 1-\alpha.
    \end{align*}
    The last inequality follows from $c_{\mathrm{max}} \geq c(x^i)$. This contradicts the minimum width of $(\hat{l}, \hat{u})$. Thus, $K \leq 1-\alpha + (c_{\mathrm{max}}/n)$.
\end{proof}}

\tb{Lemma~\ref{lem:bound} shows that when $X$ is discrete, the coverage rate of the nominal confidence band is at most $1 + \alpha + c_{\mathrm{max}}/n$. Since $\mathrm{supp}(X) = S_1 \times \cdots \times S_H$ grows exponentially with $|S_t|$ and horizon length $H$, we will have $c_{\mathrm{max}} = 1$ with a high probability, implying $c_{\mathrm{max}}/n \approx 0$ for reasonably large $|\mathrm{supp}(X)|$ and $n$. For example, suppose each $s \in \mathrm{supp}(X)$ is equally likely to be sampled. Let $V$ denote the number of possible values at each time $t$, and let $H$ be the horizon length. Then $|\mathrm{supp}(X)|$ is $V^H$, and
\begin{align*}
    P(c_{\mathrm{max}} \geq 2) = P(\exists\, i < j : x^i = x^j) 
    \leq \sum_{i < j} P(x^i = x^j) 
    = \binom{n}{2} \frac{1}{V^H},
\end{align*}
i.e., $P(c_{\mathrm{max}} \geq 2) = O(V^{-H})$. Thus, even for moderate values of $V$ and $H$, $P(c_{\mathrm{max}} = 1) \approx 1$. As a result, the under-coverage bias is present in settings where $\mathrm{supp}(X)$ is sufficiently large.
}

\section{Robust MIP}
So far, we established that an optimal solution to the nominal problem is often too narrow and thus provides a smaller coverage rate than expected. To protect against this phenomenon, we take a robust optimization perspective. Specifically, we introduce the following budget uncertainty set \cite{bertsimasPriceRobustness2004}, where the parameter $\Gamma$ takes values in $[0,1]$:
\begin{align}
    \mathcal{Z}(\Gamma) \equiv \left\{z \in [0,1]^H: \frac{1}{H}\sum_{t=1}^Hz_t \leq \Gamma \right\}
\end{align}
Intuitively, for any $z \in \mathcal{Z}(\Gamma)$, at most $\Gamma \times 100\%$ of its components can be set to 1. We then introduce the following constraints to \eqref{nominal_mip}:
\begin{align}
    \sum_{t=1}^H u_t &\geq \sum_{t=1}^H \left( q^u_t + c^u_tz_t^u \right), \quad \forall z^u \in \mathcal{Z}(\Gamma), \label{eq:robust_con_u}  \\
    \sum_{t=1}^H l_t &\leq \sum_{t=1}^H \left( q^l_t - c^l_tz_t^l \right), \quad \forall z^l \in \mathcal{Z}(\Gamma), \label{eq:robust_con_l}
\end{align}
where $c^u_t$ and $c^l_t$ are non-negative constants such that $q^u_t + c^u_t$ and $q^l_t - c^l_t$ serve as upper and lower bounds on $x^i_t$, respectively. \tb{For instance, under Assumption \ref{assn:bounds}, we may set $c^u_t = U_t - q^u_t$ and $c^l_t = q^l_t - L_t$.} These constraints help protect against under-coverage by forcing $\sum_{t=1}^H u_t$ to be larger than $\sum_{t=1}^H q^u_t$ (and $\sum_{t=1}^H l_t$ to be smaller than $\sum_{t=1}^H q^l_t$) by the extent allowed by $\Gamma$, thereby increasing the total width of the confidence band. Thus, the \textit{robust problem} is formulated as follows:
\begin{equation}\label{robust_mip} \tag{RP}
    \begin{split}
        \min_{(l,u)\in \mathcal{X}} \quad & \sum_{t=1}^H(u_t - l_t) \\
        \mbox{s.t.} \quad &\sum_{t=1}^H u_t \geq \sum_{t=1}^H \left( q^u_t + c^u_tz_t^u \right), \quad \forall z^u \in \mathcal{Z}(\Gamma), \\
        &\sum_{t=1}^H l_t \leq \sum_{t=1}^H \left( q^l_t - c^l_tz_t^l \right), \quad \forall z^l \in \mathcal{Z}(\Gamma), \\
        & \frac{1}{n}\sum_{i=1}^{n}\delta(l,u,x^i) \geq 1-\alpha,
    \end{split}
\end{equation}
where the definition of $\delta(l,u,x^i)$ can again be enforced using $n$ binary variables and an appropriate large constant. Denote by $(l^\Gamma, u^\Gamma)$ a confidence band obtained by solving the robust problem with a given $\Gamma$. We refer to this as a \textit{robust confidence band}.
\begin{rmk}
    If $\Gamma=0$, the problems \eqref{nominal_mip} and \eqref{robust_mip} are equivalent, implying $(l^0, u^0)$ would be biased towards under-coverage in the sense of Theorem \ref{thm:undercoverage}.
\end{rmk}

\subsection{Reformulation}
In this section, we present an equivalent formulation for the robust problem \eqref{robust_mip} that accounts for the infinitely many constraints introduced by \eqref{eq:robust_con_u} and \eqref{eq:robust_con_l}. 

Denote by $c^u_{(t)}$ the $t$th \textit{largest} parameter, i.e., the $t$th value from the left of the ordered parameters $c^u_{(1)} > c^u_{(2)} > \cdots > c^u_{(H)}$. By perturbing each parameter by a small amount as needed, the strict inequality holds without loss of generality. We define $c^l_{(t)}$ analogously. We will use $\lceil \cdot \rceil$ to indicate the ceiling operator, which rounds up the expression to the nearest integer. Define $t^\star \equiv \max\{\lceil \Gamma H \rceil, 1\} $ and 
\begin{align}
    \beta_{t}^u(\Gamma) &\equiv \left(c^u_{t} - c^u_{(t^\star)} \right)^+ + \Gamma c^u_{(t^\star)}, \quad \forall t,\\
    \beta_{t}^l(\Gamma) &\equiv \left(c^l_{t} - c^l_{(t^\star)} \right)^+ + \Gamma c^l_{(t^\star)}, \quad \forall t.
\end{align}
If $\Gamma=0$, then $t^\star = 1$ and we have $\beta_{t}^u(\Gamma) = \beta_{t}^l(\Gamma) = 0$ for all $t$, which will make \eqref{eq:robust_con_u2} and \eqref{eq:robust_con_l2} redundant in the reformulation. If $\Gamma=1$, then $t^\star=H$ and we have $\beta_{t}^u(\Gamma) = c^u_t$ and $\beta_{t}^l(\Gamma) = c^l_t$ for all $t$, which will widen the total width of the confidence band to the maximal extent. We present the reformulation below, where we replace constraints \eqref{eq:robust_con_u} and \eqref{eq:robust_con_l} with
\begin{align}
    \sum_{t=1}^H u_t &\geq \sum_{t=1}^H \left( q^u_t + \beta^u_t(\Gamma)  \right), \label{eq:robust_con_u2}  \\
    \sum_{t=1}^H l_t &\geq \sum_{t=1}^H \left( q^l_t - \beta^l_t(\Gamma)  \right), \label{eq:robust_con_l2}
\end{align}
respectively:
\begin{equation}\label{robust_counterpart} \tag{RP'}
    \begin{split}
        \min_{(l,u) \in \mathcal{X}} \quad & \sum_{t=1}^H(u_t - l_t) \\
        \mbox{s.t.} \quad 
        &\sum_{t=1}^H u_t \geq \sum_{t=1}^H \left( q^u_t + \beta^u_t(\Gamma)  \right), \\
        &\sum_{t=1}^H l_t \leq \sum_{t=1}^H \left( q^l_t - \beta^l_t(\Gamma)  \right), \\
        & \frac{1}{n}\sum_{i=1}^{n}\delta(l,u,x^i) \geq 1-\alpha.
    \end{split}
\end{equation}

\begin{prop} \label{prop:robust_counterpart}
    Problems \eqref{robust_mip} and \eqref{robust_counterpart} are equivalent.
\end{prop}

\begin{proof}
    We prove the equivalence by showing that constraint \eqref{eq:robust_con_u} may be replaced by \eqref{eq:robust_con_u2} without loss of optimality. The same can be shown very similarly between constraints \eqref{eq:robust_con_l} and \eqref{eq:robust_con_l2}, and we omit its proof for brevity.

    Starting from \eqref{eq:robust_con_u}, by the theorem of the alternative with appropriate primal and dual linear programs, we replace it with the following system of linear inequalities:
    \begin{align}
        \sum_{t=1}^H u_t &\geq \sum_{t=1}^H \left( q_t^u + \lambda_t^u + \Gamma \rho^u \right) \label{eq:dual_con1} \\
        \lambda_t^u + \rho^u &\geq c_t^u, \quad \forall t \label{eq:dual_con2} \\
        \lambda^u, \rho^u &\geq 0.
    \end{align}
    We further simplify these constraints by deducing the optimal values of $\lambda^u_t$ and $\rho^u$. Recall that $t^\star \equiv \max\{ \lceil \Gamma H \rceil, 1\}$ and $c^u_{(t)}$ represents the $t$th largest parameter among $c^u_t$. Suppose we enforce $\rho^u = c^u_{(t^\star)}$ as a constraint. We will show that this is not restrictive, i.e., does not harm the optimal objective value. 
    
    With $\rho^u = c^u_{(t^\star)}$, we must have at optimality:
    \begin{align*}
        \lambda^u_{(t)} = \left(c^u_{(t)} - c^u_{(t^\star)} \right)^+ = 
        \begin{cases}
            c^u_{(t)} - c^u_{(t^\star)}, &\mbox{if } t < t^\star, \\
            0, &\mbox{if } t \geq t^\star.
        \end{cases}
    \end{align*}
    Let $\epsilon > 0$ and consider the following two cases to observe their impact on the objective value:

    \textbf{Case 1}: Increase $\rho^u$ by $\epsilon$. Then for small enough $\epsilon$, we must have $\lambda^u_{(t)}=c^u_{(t)} - c^u_{(t^\star)} - \epsilon$ at optimality for all $t < t^\star$; $\lambda^u_{(t)}=0$ otherwise. The change in the objective value is
    \begin{align*}
        \Delta = \Gamma H \epsilon - (t^\star - 1)\epsilon >  0,
    \end{align*}
    since $\Gamma H > t^\star - 1$. This indicates that the new objective value is worse.

    \textbf{Case 2}: Decrease $\rho^u$ by $\epsilon$. Then for small enough $\epsilon$, we must have $\lambda^u_{(t)} = c^u_{(t)} - c^u_{(t^\star)} + \epsilon$ at optimality for $t \leq t^\star$; $\lambda^u_{(t)}=0$ otherwise. The change in the objective value is
    \begin{align*}
        \Delta = t^\star\epsilon - \Gamma H \epsilon \geq 0,
    \end{align*}
    since $t^\star \geq \Gamma H$. This indicates that the new objective value is not better.

    In either case, we cannot strictly improve the objective value. Therefore, there exists an optimal solution where $\rho^u = c^u_{(t^\star)}$ and $\lambda^u_{t} = (c^u_{t} - c^u_{(t^\star)})^+$ for all $t$. Letting $\beta^u_t(\Gamma) = \lambda^u_t + \Gamma \rho^u$ for all $t$ concludes the proof.
\end{proof}
Because \eqref{robust_mip} and \eqref{robust_counterpart} are equivalent, Proposition \ref{prop:robust_counterpart} shows that the robust problem is only slightly more complex than the nominal problem, with just two more constraints and no additional decision variables.

\subsection{Tuning $\Gamma$} \label{sec:tuning_Gamma}
By solving \eqref{robust_counterpart} for different values of $\Gamma$, we can obtain confidence bands of varying widths, starting from an under-covering solution at $\Gamma=0$ to an over-covering one at $\Gamma=1$. In this section, we present an algorithm to tune $\Gamma$. 

We cast the task of tuning $\Gamma$ as a root-finding problem. Let $f(\Gamma) = P(l^\Gamma, u^\Gamma ) - (1-\alpha)$. Then Theorem \ref{thm:undercoverage} and Assumption \ref{assn:bounds} imply $f(0) \leq 0$ and $f(1) \geq 0$, respectively. Furthermore, we note that $P(l^\Gamma, u^\Gamma)$ is continuous and increasing in $\Gamma$, which implies that there must exist $\Gamma^\star \in [0,1]$ such that $f(\Gamma^\star) = 0$. Therefore, the bisection method on $[0,1]$ will converge to $\Gamma^\star$.

In practice, the true coverage rate $ P(l^\Gamma, u^\Gamma) $ is unknown and must be estimated from finitely many sample paths. \tb{In Algorithm~\ref{alg:bisection}, we propose a bisection method with $ K $-fold cross-validation to estimate $P(l^\Gamma, u^\Gamma)$. Let $\Gamma_n^\star$ be the parameter value found using the true coverage rate $P(l^\Gamma, u^\Gamma)$, and let $\hat{\Gamma}_n$ be its estimate obtained using Algorithm \eqref{alg:bisection}. Since the algorithm accepts a fixed sample size $n$, we make the dependence on $n$ explicit in our notation. For a fixed number of folds $ K \geq 2 $, $ \hat{\Gamma}_n $ and $ \Gamma^\star_n $ become identical in the limit as the number of sample paths $ n $ and the number of iterations $ \bar{N} $ tend to infinity. This follows from the fact that any discrepancy between $ \Gamma^\star_n $ and $ \hat{\Gamma}_n $ stems solely from the difference between $ f $ and $ \hat{f} $ (defined in line 9 of Algorithm \ref{alg:bisection}). As $ n $ increases, each fold contains more samples, ensuring that the confidence bands computed from different folds, i.e., $ (l^{\hat{\Gamma}_n}_k, u^{\hat{\Gamma}_n}_k) $ in line 6, become identical in the limit. Consequently, the estimated coverage rates, i.e., $ C_k $ in line 7, also converge to the same rate for all $k$ by the law of large numbers. These conditions imply $ \hat{f} \to f $, and thus, $ \hat{\Gamma}_n \to \Gamma^\star_n $, if $ \bar{N} \to \infty $. Moreover, since $ f $ is continuous, the confidence band also becomes continuous in $ \Gamma $ in the limit. In Section~\ref{sec:case_study_VAR}, we demonstrate that even with finite samples, the algorithm yields confidence bands with accurate coverage rates.}

\begin{algorithm}
\KwIn{$n$ sample paths $x^1, \ldots, x^n$, coverage rate $1-\alpha$, number of iterations $\Bar{N}$, number of folds $K$}

\nl {\bf Initialize:} $\Gamma_a=0$, $\Gamma_b=1$, $N=0$, $K$ random partitions of $n$ sample paths, $\mathcal{P}_1, \ldots, \mathcal{P}_K$, each containing $m \equiv n/K$ paths (assume $n$ is divisible by $K$). \\
\nl {\bf while} $N < \Bar{N}$ {\bf do} \\
\nl \qquad $N = N + 1$ \\
\nl \qquad $\hat{\Gamma}_n = (\Gamma_a+\Gamma_b)/2$ \\
\nl \qquad {\bf for} $k=1,\ldots, K$ {\bf do} \\
\nl \qquad \qquad Solve \eqref{robust_counterpart} with $\hat{\Gamma}_n$ and all partitions \\
\qquad \qquad except $\mathcal{P}_k$ to obtain $(l^{\hat{\Gamma}_n}_k, u^{\hat{\Gamma}_n}_k)$. \\
\nl \qquad \qquad $C_k = \frac{1}{m}\sum_{j:x^j\in \mathcal{P}_k}\delta(l^{\hat{\Gamma}_n}_k, u^{\hat{\Gamma}_n}_k, x^j)$, i.e., coverage \\
\qquad \qquad rate over the $k$th partition $\mathcal{P}_k$ \\
\nl \qquad {\bf end for} \\
\nl \qquad $\hat{f}(\hat{\Gamma}_n) = \frac{1}{K}\sum_{k=1}^K C_k - (1-\alpha)$ \\
\nl \qquad {\bf if} $\hat{f}(\hat{\Gamma}_n) < 0$ {\bf do} \\
\nl \qquad \qquad $\Gamma_a = \hat{\Gamma}_n$ \\
\nl \qquad {\bf else} \\
\nl \qquad \qquad $\Gamma_b = \hat{\Gamma}_n$ \\
\nl \qquad {\bf end if} \\
\nl {\bf end while} \\

\KwOut{$\hat{\Gamma}_n$}
\caption{{Bisection method with $K$-fold cross-validation for estimating $\Gamma^\star_n$} \label{alg:bisection}}
\end{algorithm}

Algorithm \ref{alg:bisection} requires solving \eqref{robust_counterpart} $\Bar{N}K$ times in total to estimate $\Gamma^\star_n$. Upon estimating $\Gamma^\star_n$, we solve the robust problem once more to obtain the confidence band. In the case study in Section \ref{sec:case_study_VAR}, the total run time is only a few seconds when $n \leq 500$ at the 1\% optimality gap criterion. Nevertheless, for a potential speed-up, we observe the following.
\begin{rmk}
    The feasible set of \eqref{robust_counterpart} is decreasing in $\Gamma$.
\end{rmk}
This implies that for a fixed partition of $n$ sample paths, the confidence band obtained by solving \eqref{robust_counterpart} using $\Gamma_2$ is a feasible solution for \eqref{robust_counterpart} that uses $\Gamma_1$ for any $\Gamma_2 \geq \Gamma_1$. Therefore, solutions can be stored and used as feasible warm-starts throughout Algorithm \ref{alg:bisection} wherever applicable.

\section{Case study: estimating a confidence band for a vector autoregressive (VAR) model} \label{sec:case_study_VAR}
VAR models are commonly used in Economics and the natural sciences to examine impulse responses, which measure the impact of an external shock to one variable on others over time. To quantify the uncertainty in these responses, a confidence band must be estimated. In this section, we apply our methodology to the same two-dimensional VAR(1) example used in \cite{staszewska-bystrovaConstructingNarrowestPathwise2013} and \cite{schusslerConstructingMinimumwidthConfidence2016}. Consider 
\begin{align}
    x_t = A_0 + A_1x_{t-1} + \epsilon_t, 
\end{align}
with $x_0 = [0,0]$, where
\begin{align*}
    A_0 = 
    \begin{bmatrix}
        1 \\
        1
    \end{bmatrix}, \quad
    A_1 = 
    \begin{bmatrix}
        0.5 & 0.3 \\
        -0.6 & 1.3
    \end{bmatrix}, \quad
    \Sigma_\epsilon = 
    \begin{bmatrix}
        1 & 0.5 \\
        0.5 & 1
    \end{bmatrix}.
\end{align*}
We simulate $n$ sample paths, for different values of $n$, over 12 periods for the first variable of the VAR(1) process. We use both the baseline approach and our proposed methodology to obtain two confidence bands at $\alpha=0.1$, which we refer to as the nominal and robust confidence bands. We set the parameters $q^u_t$ and $q^l_t$ to the $(1-\alpha)$- and $\alpha$-quantiles of the $n$ sample paths at time $t$, respectively. Recall that $q^u_t + c^u_t$ represents the (practical) upper bound on the sample paths at time $t$, and thus we set it to the 1-quantile of the $n$ sample paths at time $t$. Similarly, $q^l_t - c^l_t$ represents the lower bound on the sample paths and is set to the 0-quantile of the sample paths at time $t$. Same as \cite{schusslerConstructingMinimumwidthConfidence2016}, we create four random sets $\#1, \ldots, \#4$, each with 1,000 sample paths from the VAR(1) process. In Table \ref{tab:comparison_schussler_trede}, we summarize the coverage rates by the two confidence bands.

\begin{table*}[h!]
    \centering
    \small
    \begin{tabular}{cccccc|ccccc}
    \hline
        \multirow{2}{*}{$n$} & \multicolumn{5}{c}{Nominal confidence band} & \multicolumn{5}{c}{Robust confidence band} \\
        \cline{2-11}
          & \#1 & \#2 & \#3 & \#4 & Avg. & \#1 & \#2 & \#3 & \#4 & Avg. \\
    \hline
        100 & 63.8\% & 63.1\% & 65.4\% & 65.5\% & 64.4\% & 93.1\% & 92.0\% & 91.4\% & 90.3\% & 91.7\% \\
        200 & 74.2\% & 74.6\% & 76.1\% & 77.1\% & 75.5\% & 90.1\% & 90.7\% & 90.2\% & 88.5\% & 89.9\% \\
        500 & 83.5\% & 80.9\% & 84.7\% & 83.3\% & 83.1\% & 91.0\% & 91.0\% & 89.2\% & 91.2\% & 90.6\% \\
        1,000 & 85.0\% & 83.5\% & 86.5\% & 84.6\% & 84.9\% & 90.9\% & 88.7\% & 88.9\% & 90.8\% & 89.8\% \\
        5,000 & 87.2\% & 86.9\% & 86.5\% & 89.6\% & 87.8\% & 90.4\% & 89.5\% & 89.4\% & 90.3\% & 89.9\% \\
    \hline
    \end{tabular}
    \caption{Results of the experiment for the first variable of the VAR(1) process with $\alpha = 0.1$, i.e., closer to 90\% is better. Each set $\#1, \ldots, \#4$ contains 1,000 sample paths. The table summarizes the average coverage rate over the four sets. We use $K=2$ in Algorithm \ref{alg:bisection} for $n=100$ and $n=200$ and use $K=4$ for the rest.}
    \label{tab:comparison_schussler_trede}
\end{table*}

Our approach clearly produces higher quality confidence bands whose estimated coverage rates are much closer to 90\%, especially with limited samples (e.g., $n=100$). At $n=200$, the robust confidence band already achieves an average coverage rate of approximately 90\%, whereas the nominal confidence band does not, even with 5,000 sample paths, indicating an order-of-magnitude improvement in the required sample size. Moreover, the nominal confidence band is clearly biased, as its coverage rates are always smaller than 90\%, whereas the coverage rates for the robust confidence band hover around 90\%, either slightly above or slightly below. The nominal confidence band can achieve a small enough bias with $10,000$ sample paths, as noted in \cite{schusslerConstructingMinimumwidthConfidence2016}, but the resulting MIP can take several hours to solve to optimality; even with a more relaxed optimality gap criterion, generating so many sample paths may also be very time-consuming for complex simulation models. In contrast, with a 1\% optimality gap criterion, the entire procedure of estimating $\Gamma^\star_n$ and constructing the robust solution takes roughly two and five seconds with $n=200$ and $n=500$, respectively, without using warm-starts.

\begin{figure*}[hbt!]
    \centering
    \begin{subfigure}[t]{0.4\textwidth}
        \centering
        \includegraphics[width=1\textwidth]{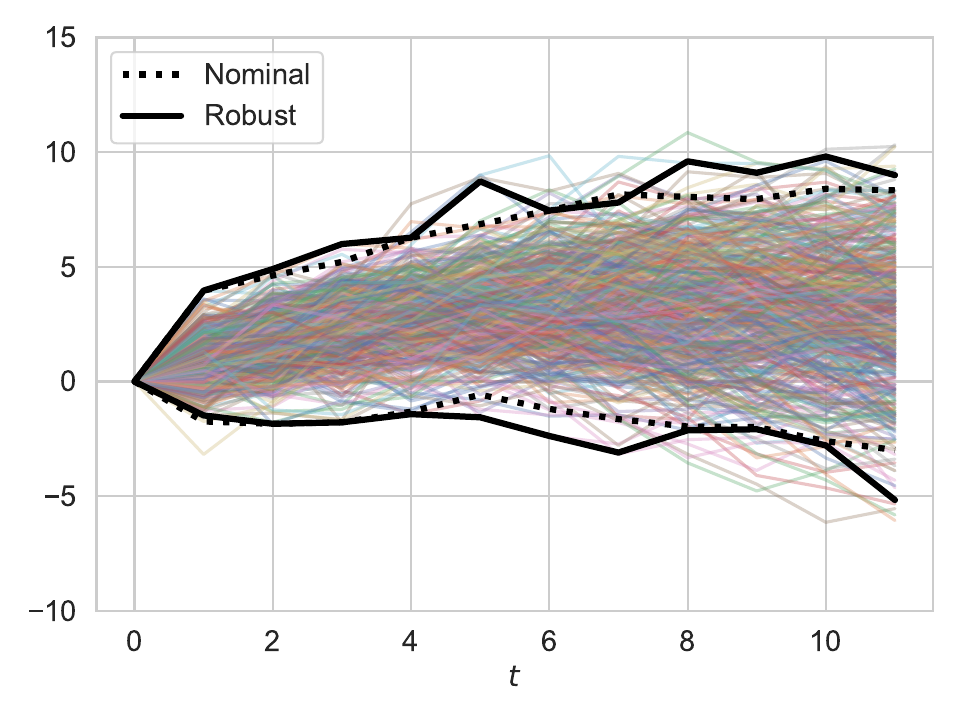}
        \caption{$n=500$}
    \end{subfigure}
    ~
    \begin{subfigure}[t]{0.4\textwidth}
        \centering
        \includegraphics[width=1\textwidth]{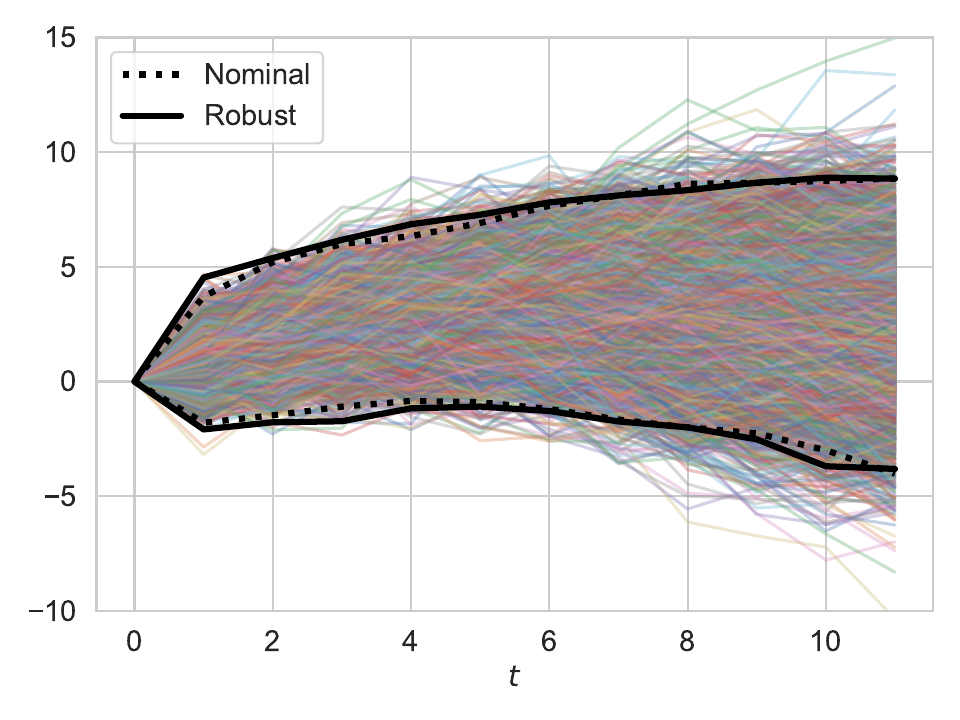}
        \caption{$n=5,000$}
    \end{subfigure}
    \caption{The nominal and robust confidence bands for the first variable of the VAR(1) process constructed using 500 (left) and 5,000 (right) sample paths.}
    \label{fig:VAR}
\end{figure*}

The differences in the nominal and robust confidence bands are illustrated in Figure \ref{fig:VAR} for $500$ and $5,000$ sample paths. With a relatively small $n$, the nominal confidence band is considerably narrower than the robust one, as it is significantly biased. As $n$ increases, the difference between the two solutions diminishes, and with $n=5,000$, they are very similar.

\section{Case study: validating a queueing model of patient flow during a Mass Casualty Event (MCE)} \label{sec:case_study_mce}
We now demonstrate how our methodology can help validate a queueing model of patient flow. We consider the Erlang-R queue in \cite{yom-tovErlangRTimeVaryingQueue2014}, which models a queue with reentrant customers who can return multiple times during their sojourn within the system. The Erlang-R queue has diverse applications with various extensions studied in the literature, see, e.g., \cite{chan2014use,chan2022dynamic}.

\begin{figure*}[hbt!]
    \centering
    \begin{subfigure}[t]{0.425\textwidth}
        \centering
        \includegraphics[width=1\textwidth]{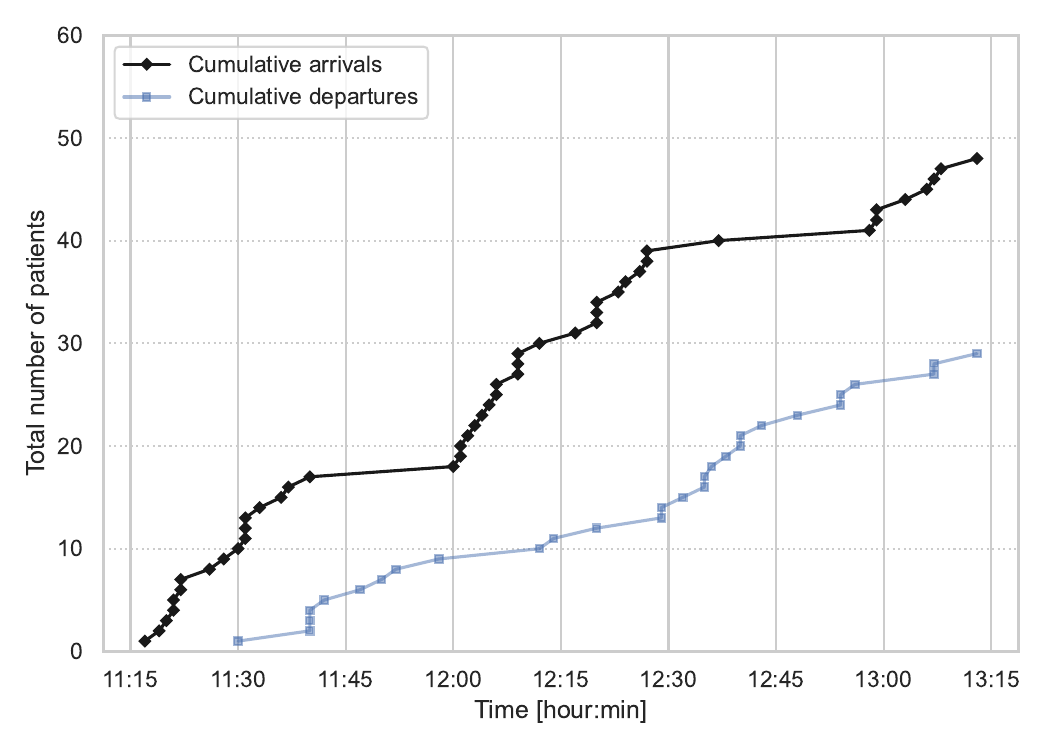}
        \caption{Actual cumulative arrivals and departures in MCE drill}
    \end{subfigure}
    ~
    \begin{subfigure}[t]{0.425\textwidth}
        \centering
        \includegraphics[width=1\textwidth]{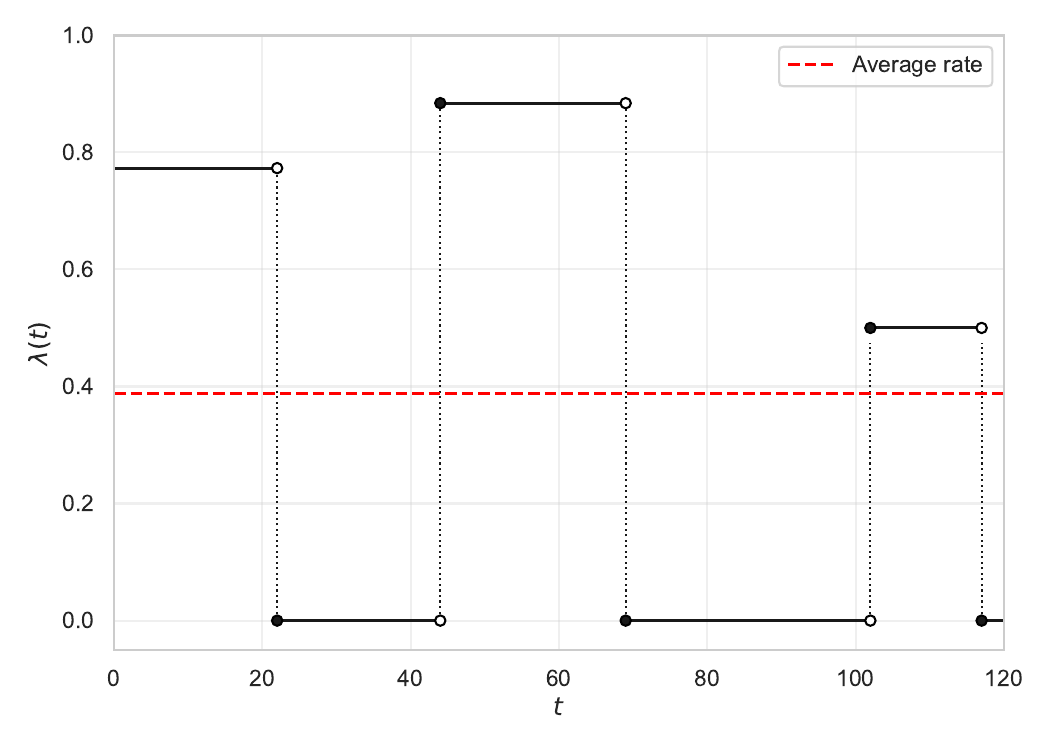}
        \caption{Estimated arrival rate function (patients per minute): $\lambda(t) = 0.773 \times 1\{0 \leq t < 22\} + 0.884 \times 1\{44 \leq t < 69\} + 0.5 \times 1\{102 \leq t < 117\}$, where $t=0$ corresponds to 11:16. The average rate is 0.388 patients per minute.}
        \label{fig:arrival_rate_function}
    \end{subfigure}
    \caption{Arrival and departure data and the estimated arrival rate function for the MCE.}
    \label{fig:arrivals_and_departures}
\end{figure*}

We consider the evolution of a hospital emergency ward during an MCE, \tb{a practical example of a discrete-valued stochastic process with large support.} We use the same data and parameters from \cite{yom-tovErlangRTimeVaryingQueue2014}, which describe a chemical MCE drill that took place in July 2010 at 11:00 and lasted until 13:15. Figure \ref{fig:arrivals_and_departures} shows the actual cumulative arrivals and departures of patients at the emergency department (left) and the estimated arrival rate function (right). Their model specified four servers. The average treatment time was 5.4 minutes ($\mu=11.06$), average time until readmission 24.6 minutes ($\delta=2.44$), and the probability of readmission $p=0.662$. \cite{yom-tovErlangRTimeVaryingQueue2014} use confidence bands from a diffusion approximation of the Erlang-R model to assess the validity of their model. Here, we use a simulation model of the Erlang-R queue to construct the confidence bands. 

\begin{figure*}[hbt]
    \centering
    \begin{subfigure}[t]{0.45\textwidth}
        \centering
        \includegraphics[width=1\textwidth]{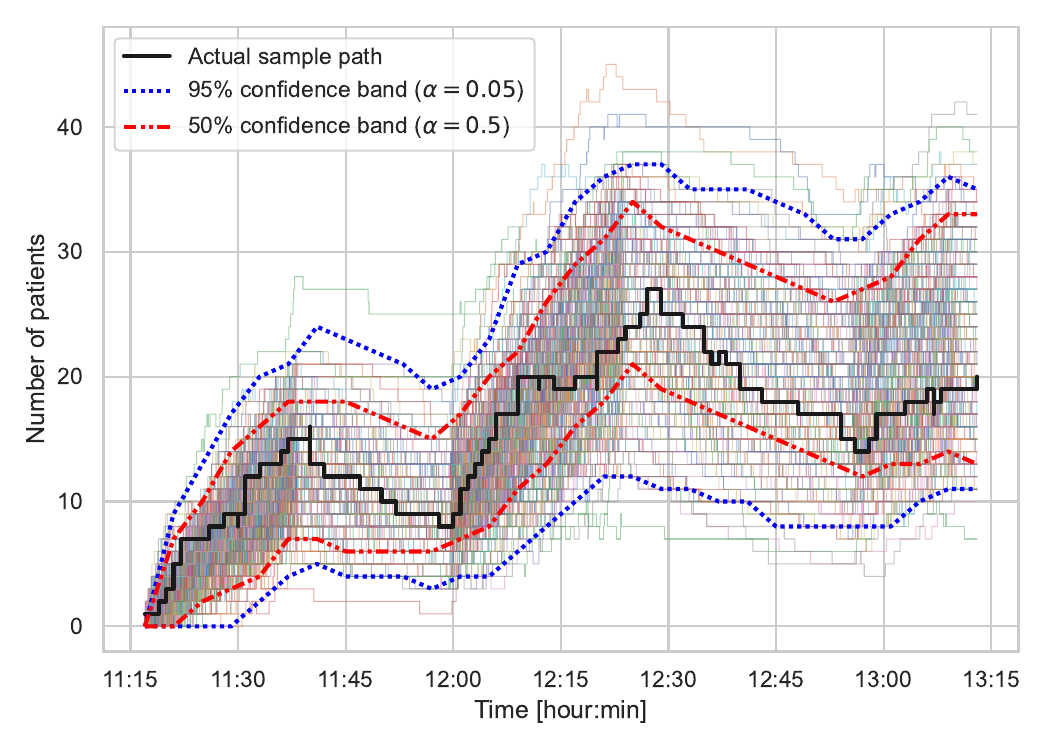}
        \caption{Assuming time-varying arrival rate function $\lambda(t)$ in Figure \ref{fig:arrival_rate_function}}
        \label{fig:erlang_R_CB}
    \end{subfigure}
    ~
    \begin{subfigure}[t]{0.45\textwidth}
        \centering
        \includegraphics[width=1\textwidth]{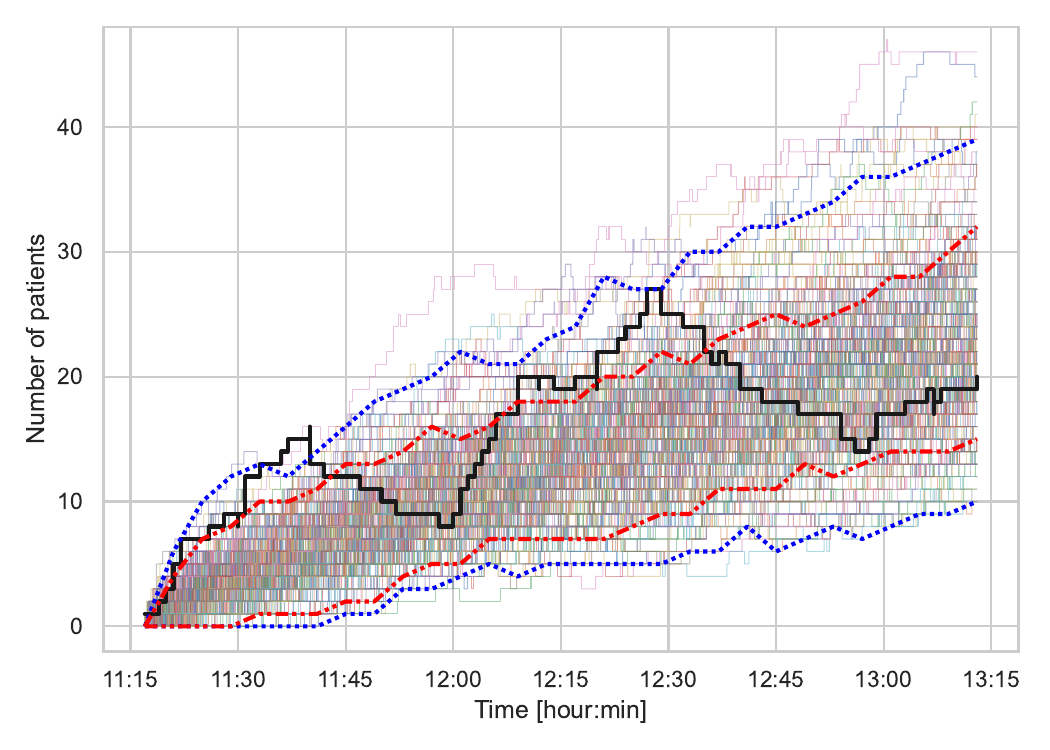}
        \caption{Assuming an avarege (stationary) arrival rate of 0.388 patients per minute}
        \label{fig:erlang_R_CB_stationary}
    \end{subfigure}
    \caption{Robust confidence band with the actual sample path.}
\end{figure*}

We generate $n=300$ sample paths. We discretize the time horizon into $H=30$ equidistant intervals and extract the value of each sample path at the start of each interval, \tb{resulting in the maximum multiplicity $c_{\mathrm{max}}=1$.} Using this data, we construct the confidence bands by running Algorithm \ref{alg:bisection} with $K=3$ and $\Bar{N}=10$ and solve \eqref{robust_counterpart} with the obtained estimate of $\Gamma^\star_n$. In Figure \ref{fig:erlang_R_CB}, we present confidence bands with $\alpha=0.5$ and $0.05$, both of which cover the actual sample path, thus supporting the validity of the specified queueing model in \cite{yom-tovErlangRTimeVaryingQueue2014}. In contrast, the diffusion approximation-based confidence bands shown in \cite{yom-tovErlangRTimeVaryingQueue2014} do not fully cover the sample path. 

Lastly, we apply our methodology to a simpler model with a stationary arrival process whose rate is the average arrival rate in Figure \ref{fig:arrival_rate_function}. The results in Figure \ref{fig:erlang_R_CB_stationary} clearly suggest that there is model mis-specification when assuming stationary arrivals, since neither band fully contains the actual sample path and they fail to capture the time-varying nature of the process.

\section{Conclusion}
We present a methodology for constructing confidence bands for discrete-time stochastic processes using a finite number of simulated sample paths. Importantly, our methodology can achieve better performance with significantly less data than existing approaches. 
Our approach improves upon existing methods by addressing optimization bias directly in the constraints, thus producing confidence bands that achieve the desired coverage rate without having to simulate many samples. We demonstrate the effectiveness of our methodology on two case studies from the literature, a vector autoregressive model and a queuing model. Note that our approach is also applicable to continuous-time processes after appropriately discretizing time. 




\bibliographystyle{elsarticle-num}
\bibliography{manuscript.bib}






\end{document}